\documentclass{amsart}
\usepackage{amsfonts,amssymb,amsmath,euscript}
\usepackage{graphicx}
\usepackage[matrix,arrow,curve]{xy}
\usepackage{wrapfig}

  \usepackage{amssymb}
\usepackage{mathrsfs}
\usepackage{euscript}
\usepackage{lipsum}

    \newtheorem{thm}{Theorem}                   
    \newtheorem{thm*}{Theorem}

    \newtheorem{lemma*}{Lemma}    

\newcommand*{\abs}[1]{\left\lvert#1\right\rvert}   

\newcommand{\hilb}{\mathcal H}

\newcommand{\clK}{\mathcal K}

\newcommand{\mbR}{\mathbb R}
\newcommand{\mbC}{\mathbb C}

\newcommand{\eps}{\varepsilon}

\begin{document}
\title{Relative stability of singular spectrum}
\author{Nurulla Azamov}
\address{Independent scholar, Adelaide, SA, Australia}
\email{azamovnurulla@gmail.com}
\begin{abstract} 
There are classical theorems of analysis which, given certain conditions on a perturbation, 
assert stability of the essential and absolutely continuous components of the spectrum of a self-adjoint operator.
Whereas the singular component is known to be highly volatile under the weakest of all 
possible perturbations, -- rank one.

In this note I announce a theorem which asserts that, nevertheless, the singular component of spectrum in an open interval is in a sense \emph{relatively stable} provided 
the limiting absorption principle (LAP) holds in the interval. One of the benefits of this result is the provision of a method for disproving LAP where it is suspected to fail.

\end{abstract}
\maketitle

\bigskip 

The limiting absorption principle (LAP) is one of the main tools in investigation 
of the spectrum of self-adjoint operators. In essence it asserts that a compact sandwiched resolvent, which sometimes is also called the Birman-Schwinger
operator, 
\begin{equation} \label{F: Tz(H0)}
    T_z(H_0) = F (H_0 - z)^{-1} F^*
\end{equation}
of a self-adjoint operator~$H_0$ converges in an appropriate topology, usually the uniform topology,
as the spectral parameter $z = \lambda+iy$ approaches the real axis for a.e. $\lambda.$ Here $F$ is some 
given or appropriately chosen closed operator which acts 
from the (complex separable) Hilbert space, $\hilb,$ where~$H_0$ acts, to possibly another Hilbert space,~$\clK.$
Once LAP is proved, the stationary approach to scattering theory allows to make conclusions about 
the absolutely continuous spectrum of the perturbed operator $H_0 + F^*JF,$ where $J$ is a bounded self-adjoint operator on~$\clK,$
see e.g. \cite{BE,YaBook,AzDaMN} and references therein.

One can consider two versions of LAP, a weak and a strong one. In a weak version the sandwiched resolvent  \eqref{F: Tz(H0)}
is assumed to have norm limits for a.e. $\lambda$ in an open interval~$I,$ while the strong LAP requires existence of uniform limit of $T_{\lambda+iy}(H_0)$ on compact subsets of~$I$
as $y \to 0^+,$ so that the resulting operator-function $T_{\lambda+i0}(H_0)$ is continuous in $\lambda\in I.$
There exists a number of methods for proving  LAP,
including the trace class method (see e.g. \cite{BE,YaBook}), a direct spectral analysis of the operator~$H_0$
where it allows one to use the so-called smooth method (see e.g. \cite{Agm,Kur}),
and most of all the celebrated conjugate operator method, initiated by E.\,Mourre, which proved to be extremely successful
in a variety of different situations, most notably for the spectral analysis of many body Schr\"odinger operators, see \cite{PSS} and \cite{AMG}.


It is known that the essential (Weyl's theorem) and a.c. (scattering theory methods) components of the spectrum
have certain stability properties under not so much restrictive assumptions on the perturbation operator.
At the same time, the singular component of the spectrum of a self-adjoint operator 
can go from one extreme, the empty set, to an everywhere dense set under a rank-one perturbation, see e.g. \cite{SiTrId2} for some striking examples. 
In this note I announce a theorem which shows that nevertheless the singular spectrum possesses a certain \emph{relative stability}
property. The premise for this property is LAP, and therefore,
this property also provides a method to disprove LAP,
by demonstrating relative instability of singular spectrum. This is important since there is a dearth of methods for disproving LAP in situations
where it is suspected to fail.

A notable case, which will be soon discussed in more detail elsewhere, where LAP is suspected to fail is the Schr\"odinger operator $H_0 = - \Delta + V_0$ with a bounded (only for simplicity, --- even this is 
apparently difficult enough) potential $V_0$ on $L^2(\mbR^d),$ and $F$ an operator of multiplication by $(1+|x|^2)^{-d/4}.$ This problem belongs to the field of potential scattering theory which has an extremely rich history going back to 1950s, see e.g. \cite{RS}, and to my knowledge it is still open. 
The fall-off $|x|^{-d}$ for the perturbing potential $V$ is well known to be special:
the cases where $F$ is an operator of multiplication by $(1+ |x|^2)^{-d/4-\eps}$ for some $\eps>0,$ or where the initial potential~$V_0$ has an additional property such as being zero or periodic, were well studied: in both cases LAP holds. 

%

\bigskip

\begin{thm}  \label{T: main thm}  Assume that 
\begin{enumerate}
  \item~$H_0$ is a self-adjoint operator on a Hilbert space $\hilb,$ and $F\colon \hilb \to\clK$ is a bounded operator with trivial kernel and co-kernel,
  \item the operator $F(H_0-\lambda-iy)^{-1}F^*$ is compact and has the limit $$F(H_0-\lambda-i0)^{-1}F^*$$
  in the uniform norm for a.e. $\lambda$ from $I,$ 
  \item $J$ is a bounded self-adjoint operator on~$\clK$  and $V = F^*JF,$ in particular,
           $V$ is also a bounded self-adjoint operator on $\hilb,$
  \item~$I$ is an open interval in $\mbR.$
\end{enumerate}
Then for any $\eps>0$ there exists a compact subset~$K$ of $I$ with $\abs{I-K}<\eps,$
such that for any $r \in \mbR$ the operator 
$$
  F E^{(s)}_K(H_r)
$$
is Hilbert-Schmidt, where $H_r = H_0+rV$ and $E^{(s)}(H)$ is the singular part of the spectral measure of~$H.$
\end{thm}

\bigskip 

Theorem \ref{T: main thm} paves the way for defining the singular spectral shift function (SSSF) $\xi^{(s)}(\lambda)$ for the pair $H_0,$ $H_0 + F^*JF$ of operators
in the interval~$I,$ see \cite{Az3v6,AzSFIES,AzSFnRI,AzDaMN}. For relatively trace class perturbations SSSF admits two other interpretations: as the singular $\mu$-invariant and as the total resonance index, see \cite{AzDaMN}. However, both the singular $\mu$-invariant
and the total resonance index require only LAP in order to be defined and in such a case they are equal, see \cite{AzDa2020}. In general, it still remains to be proved that SSSF coincides with these two other interpretations, in particular in this generality the question of integrality of SSSF is still open.

Theorem \ref{T: main thm}, a complete proof of which is still work in progress, holds for riggings~$F$ which are not necessarily bounded, but at this stage it is better not to obscure matters with more technicalities. 

\section{Sketch of proof of Theorem \ref{T: main thm}}
Here I outline some basic ideas behind the proof of Theorem \ref{T: main thm}.
Let $H_s = H_0+sV.$
Let $$R_z(H_s) = (H_s-z)^{-1}$$ be the resolvent and $$T_z(H_s) = F R_z(H_s) F^*$$ be the sandwiched resolvent. 
Recall that \emph{(coupling) resonance functions}, $$r_j(z), \ \  j = 1,2,\ldots,$$ of the pair $H_0$ and $V$ are poles of the meromorphic function 
$T_z(H_s)$ of the coupling variable~$s.$ The collection of all coupling resonance functions form an infinite-valued holomorphic function on the resolvent set of $H_0,$
see \cite{AzSFIES,AzSFIESIII,AzSFIESV} for more on this. Depending on point of view, the resonance functions may also be called resonance points. 

Let $\Theta = \{ z \in \mbC \colon z = 0 \ \text{or} \ \arg(z) \in [\pi/4, 3\pi/4 ] \}$
and for a subset $K \subset \mbR$ let $$\Theta(K) = \bigcup _{\lambda \in K} (\Theta+\lambda).$$
One essential element of the proof is the following theorem from \cite{AzSFIESIII}.

\begin{thm} \label{T0}   Assume that $H_0$ is a self-adjoint operator on a rigged Hilbert space $(\hilb,F),$ which obeys the premise of Theorem \ref{T: main thm}.  Suppose $I \subset \mbR$  is an open interval; then for any $\eps>0$ there exists a compact subset~$K$ of $I$ with $|I\setminus K| < \eps$ and~$K$ has a neighbourhood $O(K)$ in $\Theta(K)$ 
such that all coupling resonance functions $r_j(z)$ restricted to $O(K)$ either
do not take a value in $[0,1],$ or they are single-valued continuous functions without branching points in $O(K)$ (and so they have no other types of singularities too).
Moreover, the number of coupling resonances obeying the latter scenario is finite. 
\end{thm}
I will call coupling resonances from the latter scenario above \emph{impacting}. 
The set~$K$ can be chosen so that $T_z(H_0)$ admits continuous extension to $\Theta(K)$ and we will assume this. 
In fact, this is an element of the proof of Theorem \ref{T0}.

Another element of the proof is the Laurent expansion of the meromorphic function $T_z(H_s)$ in $O(K)$ which, due to Theorem \ref{T0}, has the simpler form
\begin{equation} \label{Laurent} 
   T_z(H_s) = \tilde T_z(H_s) + \sum_{j} (s - r_j(z)) ^{-1} K_z (r_j(z)),
\end{equation}
than it would otherwise, 
where the sum is over the finite number of impacting coupling resonances.
Here $K_z (r_z)$ are finite-rank operators, defined by 
\begin{equation} \label{Kz} 
   K_z(r_z) = \frac 1{2\pi i} \oint_{C(r_z)} T_z(H_s)\,ds,
\end{equation}
with $C(r_z)$ a small contour enclosing $r_z,$
 see \cite{AzSFIES},
and $\tilde T_z(H_s)$ is the part of the expansion which is holomorphic at all impacting resonance points $r_j(z).$
It can be easily shown that for impacting resonance points the functions $K_z(r_j(z))$ are continuous in $O(K)$ in the trace class norm. 

Let $C(K)$ be the set of all real-valued functions on $\mbR$ which are zero outside~$K$ and have continuous restriction to~$K.$
The third element of the proof is Stone's formula 
\begin{equation*} \label{F: phi(H)= phi*P(H) in euE}
  F\phi(H_0)F^* = \lim_{y \to 0^+} \frac 1 \pi \int \phi(\lambda) \Im T_{\lambda+iy}(H_0)\,d\lambda, \quad \phi \in C(K), 
\end{equation*}
which holds for $H_0,$ but in this form  not necessarily for $H_r, r \in (0,1],$ since $H_0$ has no singular spectrum in~$K$  while $H_r$ can. 

Further, while $T_z(H_0)$ admits continuous extension to $\Theta(K),$ the function $T_z(H_r)$ may not, if $r \in (0,1].$ 
But it can be shown that the holomorphic part $\tilde T_z(H_r)$ of the Laurent expansion \eqref{Laurent} admits continuous extension to $\Theta(K).$
Applying this to the functional 
$$
    A_r[\phi] :=   \lim_{y \to 0^+} \frac 1 \pi \int \phi(\lambda) \Im T_{\lambda+iy}(H_r)\,d\lambda,
$$
shows that it consists of two parts, the first containing $\Im \tilde T_{\lambda+iy}(H_r)$ and the second part coming from the second summand of \eqref{Laurent},
such that the first part gives an absolutely continuous operator-valued measure, and the second part, while not necessarily absolutely continuous, gives a measure 
whose singular part is trace class valued.

Finally, to complete the proof, it remains to show that $A_r[\phi]$ differs from 
$$
    F \phi(H_r)F^*
$$
by a trace class operator, which reduces to taking care of additional end-points terms in Stone's formula, see e.g. \cite{RS}.

\medskip
{\bf Remark 1.} The proof shows that the impacting coupling resonance functions are the reason the singular spectrum of $H_0$ when perturbed can encroach into the compact set $K.$
Something like this is in fact well-known in rank-one case, where there is only one coupling resonance function which in that case is Herglotz. 

\medskip
{\bf Remark 2.}  In my previous papers I used to use the operators $P_z(r_z)$ instead of $K_z(r_z)$ given in \eqref{Kz}.
The connection between these operators is that $P_z(r_z) = K_z(r_z)J.$ The operator $P_z(r_z)$ has the advantage of being an idempotent, --- it is the Riesz idempotent 
 of the eigenvalue $(s-r_z)^{-1}$ of  $T_z(H_s) J,$ which is independent from~$s.$ But using $K_z(r_z)$ allows to work with $T_z(H_s)$ instead of $T_z(H_s)J,$
 and sometimes this is more convenient. See \cite{AzDa2021} for more on this.

\medskip
{\bf Remark 3.} The reason I am announcing a result with an unfinished proof is I have no certain timeframe when it can be finished, given my current position. 
At the same time the sketch of proof provided above is detailed enough. 

\bigskip
{\it Acknowledgements.} I thank my wife for financial support during the work on this note.


\begin{thebibliography}{XXXXX}

\bibitem{Agm}   Sh.\,Agmon,   {\it Spectral properties of  Schr\"odinger operators and scattering theory,}
                    Ann. Scuola Norm. Sup. Pisa Cl. Sci. {\bf 2} (1975), 151--218

\bibitem{AMG}  W.\,O.\,Amrein, A.\,B.\,de Monvel,  V.\,Georgescu,  {\it $C_0$-groups, Commutator Methods, and Spectral Theory of $N$-Body Hamiltonians,}  Birkh\"auser, 1996

\bibitem{Az3v6}     N.\,A.\,Azamov,       {\it Absolutely continuous and singular spectral shift functions,}      
Dissertationes Math. {\bf 480} (2011), 1--102


\bibitem{AzSFIES}     N.\,A.\,Azamov,       {\it Spectral flow inside essential spectrum,}      
Dissertationes Math. {\bf 518} (2016), 1--156

\bibitem{AzSFnRI}     N.\,A.\,Azamov,       {\it Spectral flow and resonance index,}      
Dissertationes Math. {\bf 528} (2017), 1--91

\bibitem{AzSFIESIII}     N.\,A.\,Azamov,       {\it Spectral flow inside essential spectrum III: coupling resonances near essential spectrum,} arXiv: 2109.04675


\bibitem{AzSFIESV}     N.\,A.\,Azamov,       {\it Spectral flow inside essential spectrum V: on absorbing points of coupling resonances,} arXiv: 2109.14239

\bibitem{AzDaMN}     N.\,A.\,Azamov,  T.\,W.\,Daniels,  {\it Singular spectral shift function for resolvent comparable operators,} Math.\,Nachrichten  {\bf 292}  (2019), 1911--1930

\bibitem{AzDa2020}     N.\,A.\,Azamov,  T.\,W.\,Daniels,  {\it Resonance index and singular $\mu$-invariant,}  Analysis {\bf 40} (3), (2020) 151–161

\bibitem{AzDa2021}     N.\,A.\,Azamov,  T.\,W.\,Daniels,  {\it Coupling resonances and spectral properties of the product of resolvent and perturbation,}  in prep

\bibitem{BE}  M.\,Sh.\,Birman, S.\,B.\,Entina,  {\it The stationary approach in abstract scattering theory,}   Izv. Akad. Nauk SSSR, Ser. Mat.
{\bf 31}  (1967),  401--430;      English translation in Math. USSR Izv. {\bf 1} (1967)


\bibitem{Kur}  S.\,Kuroda, {\it Scattering theory for differential operators I, II,} J. Math. Soc. Japan {\bf 25} (1973),
75-104, 222-234.

\bibitem{PSS}  P.\,Perry, I.\,M.\,Sigal, B.\,Simon, {\it Spectral analysis of $N$-body Schr\"odinger operators,} Annals of Mathematics, {\bf 114} (1981), 519--567 

\bibitem{RS} M.\,Reed, B.\,Simon, {\it Methods of Modern Mathematical Physics,} vol. I--IV, Academic Press, 1972-1979

\bibitem{SiTrId2}     B.\,Simon,     {\it Trace Ideals and their Applications,} Second edition, Math. Surveys Monogr. (Amer. Math. Soc., 2005)

\bibitem{YaBook}     D.\,R.\,Yafaev,     {\it Mathematical scattering theory: general theory,}. Providence, R.\,I., AMS, 1992

\end{thebibliography}
\end{document}